\documentstyle[11pt]{article}

\newtheorem{defn}{{\bf Definition}}

\newtheorem{prop}[defn]{{\bf Proposition}}
\newtheorem{theo}{{\bf Theorem}}

\font\bbb=msbm9 scaled\magstep1

\newcommand{\RR}{\mbox{\bbb R}}

\begin{document}

\title{\bf ON POLYTOPAL UPPER BOUND SPHERES}
\author{BHASKAR BAGCHI {\normalsize AND} BASUDEB DATTA}

\date{}

\maketitle

\vspace{-7mm}

\footnotetext[1]{{\em $2010$ Mathematics Subject Classification.}
Primary 52B11, 52B05; Secondary 52B22.

\smallskip


\smallskip

The second author was partially
supported by grants from UGC Centre for Advanced Study. }

\begin{center}

\date{May 31, 2012}

\end{center}

\begin{abstract}
Generalizing a result (the case $k = 1$) due to M. A. Perles, we
show that any polytopal upper bound sphere of odd dimension $2k +
1$ belongs to the generalized Walkup class ${\cal K}_k(2k + 1)$,
i.e., all its vertex links are $k$-stacked spheres. This is
surprising since the $k$-stacked spheres minimize the face-vector
(among all polytopal spheres with given $f_0, \dots, f_{k - 1}$)
while the upper bound spheres maximize the face vector (among
spheres with a given $f_0$).

It has been conjectured that for $d\neq 2k + 1$, all $(k +
1)$-neighborly members of the class ${\cal K}_k(d)$ are tight. The
result of this paper shows that, for every $k$, the case $d = 2k
+1$ is a true exception to this conjecture.
\end{abstract}

\bigskip

We recall that a simplicial complex is said to be {\em
$l$-neighborly} if each set of $l$ vertices of the complex spans a
face. As a well known consequence of the Dehn-Sommerville
equations, any triangulated sphere of odd dimension $d = 2k + 1$
can be at most $(k + 1)$-neighborly (unless it is the boundary
complex of a simplex). A $(2k + 1)$-dimensional triangulated
sphere is said to be an {\em upper bound sphere} if it is $(k +
1)$-neighborly. This is because, by the celebrated Upper Bound
Theorem, any such sphere maximizes the face vector componentwise
among all $(2k+1)$-dimensional triangulated closed manifolds with
a given number of vertices \cite{n}.

A simplicial complex is said to be a {\em polytopal sphere} if it
is isomorphic to the boundary complex of a simplicial convex
polytope. For $n \geq 2k + 3$, the boundary complex of an
$n$-vertex $(2k + 2)$-dimensional cyclic polytope $P$ (defined as
the convex hull of any set of $n$ points on the moment curve $t
\mapsto (t, t^2, \dots, t^{2k + 2}))$ is an example of an
$n$-vertex polytopal upper bound sphere of dimension $2k + 1$.

We recall that a triangulated homology sphere $S$ is said to be
{\em $k$-stacked} if there is a triangulated homology ball $B$
bounded by $S$ all whose faces of codimension $k+1$ are in the
boundary $S$. The generalized lower bound conjecture (GLBC) due to
McMullen and Walkup \cite{mw} states that a $k$-stacked $d$-sphere
$S$ minimizes the face-vector componentwise among all triangulated
$d$-spheres $T$ such that $f_i(T) = f_i(S)$ for $0 \leq i < k$.
(Here, as usual, the {\em face-vector} $(f_0(T), \dots, f_d(T))$
of a $d$-dimensional simplicial complex $T$ is given by $f_i(T) =$
the number of $i$-dimensional faces of $T$). For polytopal spheres
$T$, this conjecture was proved by Stanley  \cite{s} and McMullen
\cite{m}. Recently, Murai and Nevo \cite{mn} proved that a
polytopal sphere (more generally, a triangulated homology sphere
with the weak Lefschetz property) satisfies equality in GLBC only
if it is $k$-stacked.

A triangulated homology ball $B$ is said to be {\em $k$-stacked}
if all its faces of codimension $k+1$ are in its boundary
$\partial B$. Thus, a triangulated (homology) $d$-sphere $S$ is
$k$-stacked if and only if there is a $k$-stacked (homology) $(d +
1)$-ball $B$ such that $\partial B = S$. As an aside, we mention
that in \cite[Theorem 2.3 (ii)]{mn}, Murai and Nevo prove:

\begin{prop}\hspace{-1.8mm}{\bf .} \label{P1}
If $S$ is a $k$-stacked triangulated homology sphere of dimension
$d\geq 2k$ then there is a unique $k$-stacked homology $(d +
1)$-ball $B$ such that $\partial B = S$. It is the largest
simplicial complex $($in the sense of set inclusion$)$ whose
$k$-skeleton agrees with that of $S$. That is, $B$ is given by the
formula \vspace{-2mm}
$$
B = \left\{\alpha \subseteq V(S) \, : \, {\alpha \choose \leq
k+1} \subseteq S\right\}.
$$
\end{prop}
(Here $V(S)$ is the vertex set of $S$ and ${\alpha \choose \leq \,
k + 1}$ denotes the set of all subsets of $\alpha$ of size $\leq k
+ 1$. Actually, Murai and Nevo give this formula with $d - k + 1$
in place of $k + 1$. But, their proof shows that it also holds
with $k + 1$ in place of $d - k + 1$, and - of course - in view of
the uniqueness statement the two formulae give the same $(d +
1)$-ball.)

This theorem gives a common generalization of Propositions 2.10,
2.11 and Corollary 3.6 of \cite{bd13v2} as well as a complete
answer to Question 6.4 of that paper. Notice that if $S$ is an
upper bound sphere of dimension $2k-1$ then $S$ is trivially
$k$-stacked. Such a sphere $S$ fails to satisfy the conclusion of
Proposition \ref{P1} unless it is the boundary complex of a
simplex. Thus, the hypothesis $d\geq 2k$ in Proposition \ref{P1}
is best possible.

We also recall that, for a $(d + 1)$-dimensional convex polytope
$P \subseteq \RR^{d + 1}$ and a point $x \not\in P$, a facet $F$
of $P$ is said to be {\em visible} from $x$ if, for any $y \in F$,
$[x, y] \cap P = \{y\}$. As a consequence of the Bruggesser-Mani
construction of a shelling order on a simplicial polytope (cf.
\cite[Theorem 8.12]{z}), we know that if $P$ is a simplicial
polytope then the facets of $P$ visible from any given point
outside $P$ form a (shellable) ball. The same holds for the facets
which are invisible from a point outside $P$.

The {\em generalized Walkup class} ${\cal K}_k(d)$ consists of the
triangulated $d$-manifolds all whose vertex-links are $k$-stacked
homology spheres. (We also note that the vertex-links of a
polytopal sphere are polytopal spheres, hence actually
triangulated spheres.) The case $k = 1$ of the following result is
due to M. A. Perles (cf. \cite[Theorem 1]{as}).

\begin{theo}\hspace{-1.8mm}{\bf .} \label{T1}
Let $S$ be a polytopal upper bound sphere of dimension $2k + 1$.
Then $S$ belongs to the generalized  Walkup class ${\cal K}_k(2k +
1)$.
\end{theo}

\noindent {\em Proof.} Let $S$ be the boundary complex of the
simplicial polytope $P$ of dimension $2k + 2$. Then $P$ is a $(k +
1)$-neighborly $(2k + 2)$-polytope. Fix a vertex $v$ of $S$, and
let $L$ be the link of $v$ in $S$. We need to prove that $L$ is
$k$-stacked. This is trivial if $P$ is a simplex. Otherwise, the
convex hull $Q$, of the vertices of $P$ excepting $v$, is again a
$(2k + 2)$-dimensional polytope. Clearly, $Q$ is also $(k +
1)$-neighborly and hence, by Radon's Theorem (cf. \cite[Page
124]{g}), $Q$ is also simplicial. Let $B$ be the pure $(2k +
1)$-dimensional simplicial complex whose facets are the facets of
$Q$ visible from $v$. By Bruggesser-Mani, $B$ is a shellable
ball.

\smallskip

\noindent {\em Claim\,:} $\partial B = L = S \cap B$.

\smallskip

Let ${\rm ast}_S(v) := \{\alpha \in S \, : \, v \not\in \alpha\}$
be the antistar of $v$ in $S$. Then ${\rm ast}_S(v)$ is a
triangulated $(2k + 1)$-ball and ${\rm ast}_S(v) \cap {\rm
star}_S(v)  = \partial({\rm ast}_S(v)) = {\rm lk}_S(v) = L$ (cf.
\cite[Lemma 4.1]{bd9}).

Let $A$ be the pure $(2k+1)$-dimensional simplicial complex whose
facets are the facets of $Q$ which are not in $B$ (i.e., invisible
from $v$). By the Bruggesser-Mani construction, $A$ is also a
shellable ball. Clearly, $\partial A =\partial B = A\cap B$.

We denote by $|A|$ the geometric carrier of $A$, i.e., the union
of the facets in $A$. If $x \in {\rm int}(|A|)$  then $x \in |S|$
and $[v, x]\cap {\rm int}(Q)$ is a non-trivial interval.
Therefore, $[v, x]\cap {\rm int}(P)$ is a non-trivial interval and
hence $x \in {\rm int}(|{\rm ast}_S(v)|)$. Thus ${\rm int}(|A|)
\subseteq |{\rm ast}_S(v)|$ and hence $|A| \subseteq |{\rm
ast}_S(v)|$. Let $y \in {\rm int}(|{\rm ast}_S(v)|)$. Let $y \in
{\rm int}(|\alpha|)$ for some $\alpha\in S$. Then $\alpha$ is a
face of $Q$. So, $y \in \partial Q$. If $y \in {\rm int}(|B|)$
then the line $l$ containing $v$ and $y$ intersect ${\rm int}(Q)$
in an interval $(y, w)$, where $y \in (v, w)$. So, $(v, w)
\subseteq {\rm int}(P)$ and $y \in (v, w)$. This is not possible
since $y \in |{\rm ast}_S(v)|$. So, $y \in |A|$. Thus, ${\rm
int}(|{\rm ast}_S(v)|) \subseteq |A|$ and hence $|{\rm ast}_S(v)|
\subseteq |A|$. So, $|{\rm ast}_S(v)| = |A|$. Since both ${\rm
ast}_S(v)$ and $A$ are subcomplexes of $S$, ${\rm ast}_S(v) =A$.
Then $\partial B = \partial A = \partial({\rm ast}_S(v)) = L$.
This proves the first equality of the claim.

Since $\partial B = L\subseteq S$, we have $\partial B\subseteq S
\cap B$. Let $y \in |S \cap B|$. If $y\in {\rm int}(|B|)$ then, by
the same argument as before, there exists $w \in P$ such that
$y\in (v, w) \subseteq {\rm int}(P)$. This is not possible since
$y \in |S|$. Therefore, $y \in |B|\setminus {\rm int}(|B|) =
|\partial B|$. Thus $|S\cap B| \subseteq |\partial B|$. Since
$\partial B \subseteq S$, this implies $S\cap B = \partial B$.
This completes the proof of the claim.

\smallskip

If $\alpha$ is a $k$-face of $B$ then $\alpha\in S$ since $S$ is
$(k+1)$-neighborly. Thus, by the claim, $\alpha \in S \cap B =
\partial B$. So, $B$ is $k$-stacked. Since, by the claim, $L =
\partial B$, it follows that $L$ is $k$-stacked. Since $L$ is the
link in $S$ of an arbitrary vertex of $S$, it follows that $S \in
{\cal K}_k(2k + 1)$. \hfill $\Box$

\medskip

In \cite{bd13v2}, we defined a {\em $k$-stellated sphere} to be a
triangulated sphere which may be obtained from the boundary
complex of a simplex by a finite sequence of bistellar moves of
index $< k$. We also defined ${\cal W}_k(d)$ as the class of all
triangulated $d$-manifolds with $k$-stellated vertex-links. An
easy induction on the number of bistellar moves used shows that
(cf. \cite[Proposition 2.9]{bd13v2})\,:

\begin{prop}\hspace{-1.8mm}{\bf .} \label{P2}
For $d \geq 2k -1$, a triangulated $d$-sphere $S$ is $k$-stellated
if and only if $S$ is the boundary of a shellable $k$-stacked $(d
+ 1)$-ball. In consequence, all $k$-stellated spheres of dimension
$\geq 2k - 1$ are $k$-stacked. Therefore, for $d \geq 2k$, ${\cal
W}_k(d) \subseteq {\cal K}_k(d)$.
\end{prop}

Thus, the proof of the above theorem shows that the polytopal
upper bound spheres of dimension $2k + 1$ are ($(k +
1)$-neighborly) members of the smaller class ${\cal W}_k(2k + 1)$.
In \cite{bd16}, we show that $(k+1)$-neighborly members of ${\cal
W}_k(d)$ are tight for $d \neq 2k + 1$. The theorem proved here
shows that the case $d = 2k + 1$ is a true exception to this
tightness criterion (since, except for the boundary complex of
simplices, no triangulated sphere can be tight).

In \cite{as}, the case $k=1$ of Theorem \ref{T1} was used to
classify the polytopal upper bound 3-spheres with 9 vertices. The
case $k=2$ of this theorem may be useful in similarly classifying
polytopal upper bound spheres of dimension 5 with few vertices.

\bigskip

{\sc Theoretical Statistics and Mathematics Unit, Indian
Statistical Institute,  Bangalore 560\,059, India}

{\em E-mail address:} {\sf bbagchi@isibang.ac.in}

\smallskip

{\sc Department of Mathematics, Indian Institute of Science,
Bangalore 560\,012,  India}

{\em E-mail address:} {\sf   dattab@math.iisc.ernet.in}


{\footnotesize

\begin{thebibliography}{99}
\bibitem{as}
A. Altshuler and L. Steinberg, Neighborly 4-polytopes with 9
vertices, {\em J. Combin. Theory} (A) {\bf 15} (1973), 270--287.

\bibitem{bd9}
B. Bagchi, B. Datta, Lower bound theorem for normal
pseudomanifolds, {\em  Expositiones Math.} {\bf 26} (2008),
327--351.

\bibitem{bd13v2}
B. Bagchi and B. Datta, On stellated spheres, shellable balls,
lower bounds and a combinatorial criterion for tightness,
arXiv:\,1102.0856\,v2, 2011, 46 pages.

\bibitem{bd16}
B. Bagchi and B. Datta, On stellated spheres and a tightness
criterion for combinatorial manifolds (preprint).

\bibitem{g}
B. Gr\"{u}nbaum, {\em Convex Polytopes} - 2nd ed. (GTM 221),
Springer-Verlag, New York, 2003.

\bibitem{m}
P. McMullen, On simple polytopes, {\em Invent. Math.} {\bf 113}
(1993), 419--444.

\bibitem{mw} P. McMullen, and D. W. Walkup, A generalized
lower-bound conjecture for simplicial polytopes, {\em Mathematica}
{\bf 18} (1971), 264--273.

\bibitem{mn}
S. Murai, and E. Nevo, On the generalized lower bound conjecture
for polytopes and spheres arXiv:\,1203.1720\,v2, 2012, 14 pages.

\bibitem{n}
I. Novik, Upper bound theorems for homology manifolds, {\em Israel
J. Math.} {\bf 108} (1998), 45--82.

\bibitem{s}
R. P. Stanley, The number of faces of a simplicial convex
polytope, {\em Advances in Math} {\bf 11} (1980), 236--238.

\bibitem{z} G. M. Ziegler, {\em Lectures on Polytopes},
Springer-Verlag, New York, 1995.

\end{thebibliography}
}
\end{document}